\documentclass[12pt,a4paper]{amsart}
\usepackage[T2A]{fontenc}
\usepackage[cp1251]{inputenc}
\usepackage[english]{babel}
\usepackage{amssymb, amsmath, latexsym}
\usepackage{mathtext,cite,enumerate,float,amsthm}
\usepackage{graphicx}
\renewcommand{\le}{\leqslant}
\renewcommand{\ge}{\geqslant}
\newcommand{\complex}{\mathbb{C}}
\newcommand{\diam}{\text{\,\rm diam\,}}
\newcommand{\dist}{\text{\,\rm dist\,}}
\renewcommand{\Re}{\text{\rm Re\,}}
\renewcommand{\Im}{\text{\rm Im\,}}
\newcommand{\spt}{\text{\rm spt\,}}

\newtheoremstyle{citedth}%
  {5pt}
  {5pt}
  {\itshape}
  {}
  {\bfseries}
  {.}
  {.3em}
  {\thmname{#1} \thmnumber{#2} \thmnote{\normalfont#3}}

\newcommand*{\medcup}{\mathbin{\scalebox{1.2}{\ensuremath{\bigcup}}}}%

\theoremstyle{citedth}
\newtheorem{theoremA}{Theorem}

\theoremstyle{theorem}
\newenvironment{Proof}
{\noindent{\bf Proof.}} {\hfill$\scriptstyle\blacksquare$ \medskip}
\newtheorem{theorem}{Theorem}

\newtheorem*{theorem*}{Theorem~\ref{new_theorem2}}
\newtheorem{lemma}{Lemma}

\theoremstyle{definition}

\newtheorem{remark}{Remark}

\usepackage{geometry} 
\begin{document}
\title[Singular integrals unsuitable for the curvature method]
{Singular integrals \\ unsuitable for the curvature method \\
whose $L^2$-boundedness still implies rectifiability}
\author{Petr Chunaev}
\address{Departament de Matem\`{а}tiques, Universitat Aut\`{o}noma de Barcelona, Edifici C
08193 Bellaterra (Barcelona) Catalonia} \email{chunayev@mail.ru,
chunaev@mat.uab.cat}
\author{Joan Mateu}
\address{Departament de Matem\`{а}tiques, Universitat Aut\`{o}noma de Barcelona, Edifici C
08193 Bellaterra (Barcelona) Catalonia} \email{mateu@mat.uab.cat}
\author{Xavier Tolsa}
\address{ICREA and Departament de Matem\`{а}tiques, Universitat Aut\`{o}noma de Barcelona, Edifici C
08193 Bellaterra (Barcelona) Catalonia} \email{xtolsa@mat.uab.cat}
\thanks{The research was supported by the ERC grant 320501 of the European
Research Council (FP7/2007-2013).}

\begin{abstract}
The well-known curvature method initiated in works of Melnikov and
Verdera is now commonly used to relate the $L^2(\mu)$-boundedness of
certain singular integral operators to the geometric properties of the support of
measure~$\mu$, e.g. rectifiability. It can be applied however only if
Menger curvature-like permutations, directly associated with the
kernel of the operator, are non-negative. We give an example of an
operator in the plane whose corresponding permutations change sign
but the $L^2(\mu)$-boundedness of the operator still implies that
the support of~$\mu$ is rectifiable. To the best of our knowledge, it is the
first example of this type. We also obtain several related results with
Ahlfors-David regularity conditions.
\end{abstract}

\maketitle

\section{Introduction}
We start with necessary notation and background facts. Note that we
work only in the complex plane and therefore usually skip dimension
markers in definitions.

Let $E\subset \mathbb{C}$ be a Borel set and $B(z,r)$ be an open disc with center $z\in \mathbb{C}$ and radius~${r>0}$.
 We denote by $\mathcal{H}^1(E)$ the $1$-dimensional Hausdorff measure of $E$. A set $E$ is called \textit{rectifiable} if it is contained in a countable union of Lipschitz graphs, up to a set of $\mathcal{H}^1$-measure zero. A set $E$ with $\mathcal{H}^1(E)<\infty$ is called \textit{purely unrectifiable} if it intersects any Lipschitz graph in a set of $\mathcal{H}^1$-measure zero.

We say that $K(\cdot,\cdot):\mathbb{C}^2\setminus\{(z,\zeta)\in \mathbb{C}^2: z=\zeta\}\to \mathbb{C}$
 is a Calder\'{o}n-Zygmund kernel if there exist constants $C>0$ and $\eta\in(0,1]$
 such that for all $z,\zeta \in \mathbb{C}$, $z\neq \zeta$, it holds
 that  $|K(z,\zeta)|\le C|z-\zeta|^{-1}$ and
$$
|K(z,\zeta)-K(z',\zeta)|,|K(\zeta,z)-K(\zeta,z')|\le
C\frac{|z-z'|^\eta}{|z-\zeta|^{1+\eta}}\quad \text{if}\quad
|z-z'|\le \tfrac{1}{2}|z-\zeta|.
$$

Given a fixed positive Radon measure $\mu$ on $\mathbb{C}$, a
Calder\'{o}n-Zygmund kernel $K$  and an ${f\in L^1(\mu)}$,
 we define the truncated \textit{Calder\'{o}n-Zygmund operator} (\textit{CZO}) as
\begin{equation}
\label{integral_operator} T_{K,\varepsilon}f(z):=\int_{E \setminus
B(z,\varepsilon)}
f(\zeta)K(z-\zeta)d\mu(\zeta),\qquad E=\spt \mu, \qquad \varepsilon>0.
\end{equation}
We do not define the CZO $T_K$ explicitly because several delicate
problems, such as the existence of the principal value, might arise.
On the contrary, the integral in (\ref{integral_operator})  always
converges absolutely and thus the principal value problem can be
avoided. Nevertheless, we  say that $T_K$ is
\textit{$L^2(\mu)$-bounded} if the operators $T_{K,\varepsilon}$ are
$L^2(\mu)$-bounded uniformly on $\varepsilon$.

\textit{How to relate the $L^2(\mu)$-boundedness of a certain CZO to the geometric properties of the support of $\mu$} is an old problem in harmonic analysis. It stems from Calder\'{o}n's paper~\cite{Calderon} where it is proved that the Cauchy transform,
i.e. the CZO with $K(z)=1/z$, is $L^2(\mathcal{H}^1\lfloor E)$-bounded if $E$ is a Lipschitz graphs
with small slope. Later on, Coifman, McIntosh and Meyer \cite{CMM}
removed the small Lipschitz constant assumption. In \cite{David}
David fully characterized rectifiable \textit{curves} $\Gamma$, for which the
Cauchy transform is $L^2(\mathcal{H}^1\lfloor \Gamma)$-bounded. 
These results led to further development of tools for understanding
the above-mentioned problem.

A new quantitative characterization of rectifiability in terms of the so-called $\beta$-numbers introduced
by Jones~\cite{Jones} and the concept of uniform rectifiability
proposed by David and Semmes~\cite{David_Semmes,DavSem_Ast} are
among these tools. Several related definitions for the plane are in
order. (We refer the reader to~\cite{David_Semmes,DavSem_Ast} for
definitions and results in the multidimensional case). A Radon
measure $\mu$ on $\mathbb{C}$ is called (\textit{$1$-dimensional}) \textit{Ahlfors-David regular} (or \textit{AD-regular}, for short) if it satisfies the inequalities
\begin{equation*}
\label{AD-regular}
C^{-1}r\le \mu(B(z,r))\le Cr,\qquad\text{where}\qquad z\in \spt \mu,\qquad 0<r<\diam (\spt \mu)
\end{equation*}
and $C>0$ is some fixed constant. A measure $\mu$ is called \textit{uniformly rectifiable} if it is AD-regular and $\spt \mu$ is contained in an AD-regular curve.

The well-known \textit{David-Semmes problem} is stated in the plane
as follows: \textit{does the $L^2(\mu)$-boundedness of  the Cauchy transform is sufficient for the uniform
rectifiability of the AD-regular measure $\mu$?}  This problem  was settled by Mattila, Melnikov and Verdera in
\cite{MMV}:
\begin{theoremA}[\cite{MMV}]
\label{theorem_AD_MMV} Let $\mu$ be an AD-regular measure on
$\mathbb{C}$. The measure $\mu$ is uniformly rectifiable if and only
if the Cauchy transform is $L^2(\mu)$-bounded.
\end{theoremA}

Note that an analogous problem in higher dimensions in the codimension $1$ was more recently solved
by Nazarov, Tolsa and Volberg in~\cite{NTV}.

The proof of Theorem~\ref{theorem_AD_MMV} relied on the so-called \textit{curvature
method} that was new at that time but soon became very influential
in solving many long-standing problems (see
\cite{Tolsa_book} and especially historical remarks there). Let us describe the
heart of the method. Given pairwise distinct points $z_1,z_2,z_3\in
\mathbb{C}$, their \textit{Menger curvature} is
$$
c(z_1,z_2,z_3)=\frac{1}{R(z_1,z_2,z_3)},
$$
where $R(z_1,z_2,z_3)$ is the radius of the circle passing through
$z_1$, $z_2$ and $z_3$
 (with $R(z_1,z_2,z_3)=\infty$ and $c(z_1,z_2,z_3)=0$ if the points are collinear).
This geometric characteristic is closely related to the Cauchy kernel
as shown by Melnikov \cite{M}:
\begin{equation}
\label{c^2}
 c(z_1,z_2,z_3)^2=\sum_{s\in
\mathfrak{S}_3}\frac{1}{(z_{s_2}-z_{s_1})\overline{(z_{s_3}-z_{s_1})}},
\end{equation}
where $\mathfrak{S}_3$ is the group of permutations of three
elements. Moreover, Melnikov also introduced a notation of
\textit{the curvature of a Borel measure} $\mu$:
\begin{equation}
\label{c2(mu)}
 c^2(\mu)=\iiint
c(z_1,z_2,z_3)^2\;d\mu(z_1)\,d\mu(z_2)\,d\mu(z_3).
\end{equation}
One can consider $c^2_{\varepsilon}(\mu)$, a truncated version of
$c^2(\mu)$, which is defined via the triple integral in (\ref{c2(mu)}) but over the
set
$$
\{(z_1,z_2,z_3)\in \mathbb{C}^3: |z_k-z_j|\ge \varepsilon>0, \quad
1\le k,j\le 3, \quad j\neq k\}.
$$
If $\mu$ in addition has linear growth, i.e. $\mu(B(z,r))\le Cr$ for
all $z\in \spt \mu$, then the relation between the curvature and the
$L^2(\mu)$-norm of the Cauchy transform (of measure) is specified by
the following identity due to Melnikov and Verdera \cite{MV}:
\begin{equation}
\label{Melnikov-Verdera} \int\left|\int_{\mathbb{C}\setminus
B(z,\varepsilon)}\frac{d\mu(\zeta)}{\zeta-z}\right|^2d\mu(z)=
\tfrac{1}{6}c_\varepsilon^2(\mu)+\mathcal{R}_\varepsilon(\mu),\qquad
|\mathcal{R}_\varepsilon(\mu)|\le C\mu(\mathbb{C}).
\end{equation}

The formulas (\ref{c^2}) and (\ref{Melnikov-Verdera}), generating
the curvature method, are remarkable in the sense that they relate
an analytic notion (the Cauchy transform) with a metric-geometric
one (the curvature). It is however very important here that  the
permutations in (\ref{c^2}) are always non-negative.

Later on, Theorem~\ref{theorem_AD_MMV} was pushed even further by
David and L\'{e}ger \cite{L,David_revista}. They essentially used
the non-negativity of (\ref{c^2}) in the proof of the following
assertion.
\begin{theoremA}[\cite{L}]
\label{theorem_Leger} Let $E\subset\mathbb{C}$ be a Borel set such
that $0<\mathcal{H}^1(E)<\infty$.  If the Cauchy transform is
$L^2(\mathcal{H}^1\lfloor E)$-bounded, then $E$ is rectifiable.
\end{theoremA}
Note that the  $L^2(\mathcal{H}^1\lfloor E)$-boundedness of the
Cauchy transform and the identity~(\ref{Melnikov-Verdera}) imply
that $c^2(\mathcal{H}^1\lfloor E)<\infty$. Consequently, to prove
Theorem~\ref{theorem_Leger} it is enough to show that
$c^2(\mathcal{H}^1\lfloor E)<\infty$ and this was actually done in
\cite{L}.

Until recently, very few things were known in this direction beyond
the CZO associated to the Cauchy kernel and its coordinate
parts $\Re z/|z|^2$ and $\Im z/|z|^2$,  see \cite{MMV,CMPT}. But recently
Chousionis, Mateu, Prat and Tolsa \cite{CMPT} (see also \cite{CMPTcap}) extended
Theorems~\ref{theorem_AD_MMV} and~\ref{theorem_Leger} to the CZOs
associated with the kernels
\begin{equation}
\label{kappa_n} \kappa_n(z):=\frac{(\Re z)^{2n-1}}{|z|^{2n}}, \qquad
n\in \mathbb{N},
\end{equation}
thus providing for $n\ge 2$ the first non-trivial example of CZOs
with the above-mentioned properties but not directly related to the
Cauchy transform (for $n=1$ one gets $\Re z/|z|^2=\Re (1/z)$). Note
that the results in \cite{CMPT} require a bit different notation
than in Theorems~\ref{theorem_AD_MMV} and~\ref{theorem_Leger}.
Namely, given a \textit{real-valued} Calder\'{o}n-Zygmund kernel
$K$, one has to consider the following permutations that substitute
the curvature (\ref{c^2}):
\begin{equation}
\label{permutation_Main}
\begin{split}
p_K(z_1,z_2,z_3):=&K(z_{1}-z_{2})K(z_{1}-z_{3})+K(z_{2}-z_{1})K(z_{2}-z_{3})\\&+K(z_{3}-z_{1})K(z_{3}-z_{2}).
\end{split}
\end{equation}
Analogously to (\ref{c2(mu)}), for any Borel measure $\mu$ set
\begin{equation}
\label{p_K(mu)} p_K(\mu)=\iiint
p_K(z_1,z_2,z_3)\;d\mu(z_1)\,d\mu(z_2)\,d\mu(z_3).
\end{equation}
One can also define $p_{\varepsilon,K}(\mu)$, the truncated version
of $p_K(\mu)$, in an obvious way.

In the case of kernels (\ref{kappa_n}) as in \cite{CMPT} one puts
$K(z)=\kappa_n(z)$ in (\ref{permutation_Main}) and (\ref{p_K(mu)}).
It is shown in \cite{CMPT} that the permutations
$p_{\kappa_n}(z_1,z_2,z_3)$ are non-negative for all triples
$(z_1,z_3,z_3)\in \mathbb{C}^3$ and this is appreciably used in a
curvature-like method in \cite{CMPT}.

In \cite{Chunaev2016}, kernels of the form
\begin{equation}
\label{kernels_n_N} K_t(z):=\kappa_N(z)+t\cdot \kappa_n(z),\qquad
n\le N, \qquad n,N\in \mathbb{N},\qquad t\in \mathbb{R},
\end{equation}
i.e. linear combinations of the kernels (\ref{kappa_n}) of different
order, were introduced. Clearly, one obtains a kernel of the form
(\ref{kappa_n}) from (\ref{kernels_n_N}) when $n=N$ (and $t\neq -1$)
or $t=0$. It turns out that this slight modification of the kernel
leads to a diverse behaviour of the corresponding CZO depending on
the parameter $t$. For example, it is shown in \cite{Chunaev2016}
that if $t$ belongs to the set
\begin{equation}
\label{Omega} \Omega(n,N):=\left\{\begin{array}{l}
\{0\}\cup\mathbb{R}\setminus \left(-\tfrac{1}{2}\left(3+\sqrt{9-4\tfrac{N}{n}}\right); 2-\tfrac{N}{n}\right) \text{ if } n<N\le 2n, \\
\{0\}\cup\mathbb{R}\setminus
\left(-\frac{1}{2}\left(\sigma_{n,M}+\sqrt{\sigma_{n,M}^2-4\frac{N}{n}}\right);
\sigma_{n,M}-3\right) \text{ if } N\ge 2n,
\end{array}
\right.
\end{equation}
where $\sigma_{n,M}:=3+(\tfrac{N}{n}-2)\sqrt{N-2n}$, then
\begin{equation}
\label{new_theorem1}
p_{K_t}(z_1,z_2,z_3)\ge 0 \quad \text{for all} \quad
(z_1,z_2,z_3)\in \mathbb{C}^3.
\end{equation}
Moreover, taking into account this property and using a curvature-like method, the following Theorem~\ref{theorem_Leger} type result is proved in \cite{Chunaev2016}.
\begin{theoremA}[\cite{Chunaev2016}]
\label{new_theorem2_CHU} Let $E\subset \mathbb{C}$ be a Borel set
such that $0<\mathcal{H}^1(E)<\infty$. If the~CZO $T_{K_t}$ with
$t\in \Omega(n,N)$ is $L^2(\mathcal{H}^1\lfloor E)$-bounded,
then
$E$ is rectifiable.
\end{theoremA}
From the other side, it is shown in \cite{Chunaev2016} that there
exist triples $(z_1,z_2,z_3)$ such that $p_{K_t}(z_1,z_2,z_3)$
change sign if $t$ belongs to the interval
 \begin{equation}
\label{omega}
 \omega(n,N):=(-N/n;0)
\end{equation}
Obviously, $\omega(n,N)\subseteq \mathbb{R}\setminus \Omega(n,N)$.
Note also that $\omega(n,2n)=(-2;0)=\mathbb{R}\setminus
\Omega(n,2n)$. For this reason, a curvature-like method cannot be
applied directly for $t\in \omega(n,N)$. Moreover, it follows from
 Huovinen's result in \cite{H} that Theorem~\ref{new_theorem2_CHU}
fails for ${t=-1\in \omega(n,N)}$ in the sense that there exists an
\textit{AD-regular} purely unrectifiable set $E$ with
$0<\mathcal{H}^1(E)<\infty$ such that the operator $T_{K_t}$ with
$t=-1$  is $L^2(\mathcal{H}^1\lfloor E)$-bounded. In the case
$(n,N)=(1,2)$, i.e. for the kernels
\begin{equation}
\label{kernels_1_2}
k_t(z):=\frac{(\Re z)^{3}}{|z|^{4}}+t\cdot\frac{\Re z}{|z|^{2}},\qquad
t\in \mathbb{R},
\end{equation}
even more is known due to Jaye and Nazarov \cite{JN}. Namely,  for
$t=-3/4\in \omega(1,2)$  there also exists a purely unrectifiable
(but not AD-regular) set $E$ such that $T_{k_t}$ is
$L^2(\mathcal{H}^1\lfloor E)$-bounded. For the details see also
\cite[Remark 2]{Chunaev2016}.

Thus we come to the question of what happens when $t\in\omega(n,N)$,
i.e. the permutations $p_{K_t}(z_1,z_2,z_3)$ change sign, and
curvature-like methods as in \cite{MMV,L,CMPT,Chunaev2016} do not
work. In this paper a partial answer is given in the case of
kernels~(\ref{kernels_1_2}). Namely, we show that for
$t\in(-2;-\sqrt{2})\subset\omega(1,2)$ the analogues of
Theorems~\ref{theorem_AD_MMV} and~\ref{theorem_Leger} are still
valid (a plausible conjecture for the kernels~(\ref{kernels_n_N})
with $t\in\omega(n,N)$ is also stated). To the best of our
knowledge, this is the first example of kernels with this property
in the plane. We also establish an analogue of
Theorem~\ref{theorem_AD_MMV} for the kernels~(\ref{kernels_n_N})
with $t\in\Omega(n,N)$. The corresponding results are given in the
next section. 


\section{Main results}

The following two theorems are analogues of
Theorems~\ref{theorem_AD_MMV} and~\ref{theorem_Leger} for the
kernels~(\ref{kernels_1_2}) with $t\in(-\sqrt{2};-2)$, whose corresponding
permutations change sign and a~curvature-like method cannot be
applied directly.
We will prove them in Section~\ref{section3} by exploiting sharp
estimates for permutations related to the kernels~(\ref{kappa_n})
but not to the ones in~(\ref{kernels_1_2}). Recall that
$\omega(1,2)=(-2;0)$, see (\ref{omega}).
\begin{theorem}
\label{theorem_AD_particular} Let $\mu$ be an AD-regular measure on
$\mathbb{C}$ and $T_{k_t}$ the CZO associated with the
kernel~$(\ref{kernels_1_2})$, where $t\in(-2;-\sqrt{2})\subset
\omega(1,2)$. The measure $\mu$ is uniformly rectifiable if and only
if $T_{k_t}$ is $L^2(\mu)$-bounded.
\end{theorem}
Note that this theorem fails if $t=-1\in\omega(1,2)$. It follows
from the aforementioned Huovinen's result~\cite{H}.
\begin{theorem}
\label{theorem_Leger_particular} Let $E\subset \mathbb{C}$ be a
Borel set such that $0<\mathcal{H}^1(E)<\infty$, and $T_{k_t}$ the
CZO associated with the kernel~$(\ref{kernels_1_2})$,
 where $t\in(-2;-\sqrt{2})\subset \omega(1,2)$. If $T_{k_t}$  is $L^2(\mathcal{H}^1\lfloor
 E)$-bounded,
then $E$ is rectifiable.
\end{theorem}

\begin{figure}
  \includegraphics[width=15cm]{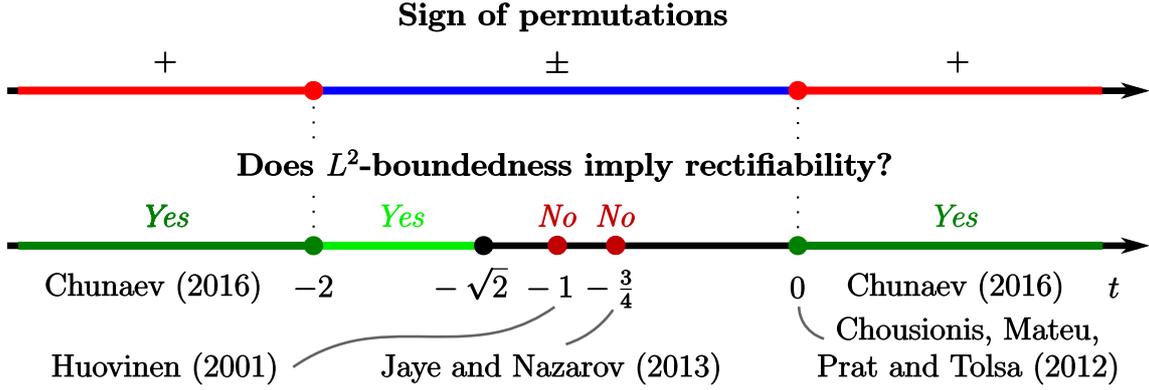}\\
 \caption{The overall picture for CZO associated with the
 kernels~$k_t$}
  \label{FIG1}
\end{figure}

This theorem supplements Theorem~\ref{theorem_Leger} type
results about CZO associated with the kernels~$k_t$ (see Figure~\ref{FIG1}). By
\cite{Chunaev2016}, if $t\notin (-2;0)$, then the permutations $p_{k_t}$ are non-negative and the
$L^2(\mathcal{H}^1\lfloor E)$-boundedness of $T_{k_t}$ implies that
$E$ is rectifiable by a curvature-like method. According to
\cite{Chunaev2016}, the permutations $p_{k_t}$ for $t\in
(-2;0)$ change sign, and by \cite{H,JN}
if $t=-1$ or $t=-3/4$, then the
operator $T_{k_t}$ is $L^2(\mathcal{H}^1\lfloor E)$-bounded but $E$
is not rectifiable. The interval $(-2;-\sqrt{2})$ corresponds to
Theorem~\ref{theorem_Leger_particular} of this paper.
\begin{remark}
As we will see at the end of Section~\ref{section3}, it is plausible
 that analogues of Theorems~\ref{theorem_AD_particular}
and~\ref{theorem_Leger_particular} are valid for the kernels~(\ref{kernels_n_N}) with $|t|>\sqrt{N/n}$. Note that in particular
$(-N/n;-\sqrt{N/n})\subset \omega(n,N)$, i.e. for $t$ from this
interval the corresponding permutations change sign. 
\end{remark}

We now formulate a Theorem~\ref{theorem_AD_MMV} type result for the
kernels (\ref{kernels_n_N}) in the case where $t\in \Omega(n,N)$,
i.e. the corresponding permutations are non-negative  (see
(\ref{Omega}) and (\ref{new_theorem1})).
\begin{theorem}
\label{theorem_AD_general} Let $\mu$ be an AD-regular measure on
$\mathbb{C}$ and $T_{K_t}$ the CZO  associated with the
kernel~$(\ref{kernels_n_N})$, where $t\in \Omega(n,N)$. The measure
$\mu$ is uniformly rectifiable if and only if $T_{K_t}$ is
$L^2(\mu)$-bounded.
\end{theorem}

Since the permutations are non-negative here, we can use a
curvature-like method. The proof that will be given in
Section~\ref{section_AD_general} is more or less analogous to the one used for the kernels~(\ref{kappa_n}) in \cite[Section~8]{CMPT}.

\bigskip

Let us say a few words about the notation in this paper. As usually, $C$ stands for a positive constant which may change its value in different occurrences. Sometimes $C$ may depend on some parameters and then we indicate it by writing, for instance, $C(\varepsilon)$ or $C_\varepsilon$, where $\varepsilon$ is a parameter. On the other hand, constants with subscripts, such as $\eta_1$ or $\theta_0$, retain their values at different places. The notation $A\lesssim B$ for positive $A$ and $B$ means that there is a positive constant $C$ such that $A\le C\,B$. If this $C$ depends on a parameter, say, $\varepsilon$, we write $A\lesssim_{\varepsilon} B$. Also, $A\approx B$ is equivalent to $A\lesssim B \lesssim A$.

\section{The proof of Theorems~\ref{theorem_AD_particular} and~\ref{theorem_Leger_particular}}
\label{section3}

Recall that
$$
\kappa_1(z)=\frac{\Re z}{|z|^{2}},\qquad \kappa_2(z)=\frac{(\Re
z)^{3}}{|z|^{4}}\qquad \text{and} \qquad k_t(z)=\kappa_2(z)+t\cdot
\kappa_1(z).
$$
The following result from \cite{CMPT} will be necessary below.
\begin{lemma}[Proof of Proposition 2.1 in \cite{CMPT}]
\label{lemma_CMPT} Given $u=(x,y)$ and $v=(a,b)$ in $\mathbb{C}$,
\begin{equation}
\label{representation_CMPT} p_{\kappa_m}(0,u,v)=\sum_{k=1}^m
\binom{m}{k}
\frac{(ax(x-a))^{2(m-k)}}{|u|^{2m}|v|^{2m}|u-v|^{2m}}\,h_k(u,v),
\end{equation}
where
$h_k(u,v):=(ax)^{2k-1}(y-b)^{2k}+(x(x-a))^{2k-1}b^{2k}+(a(a-x))^{2k-1}y^{2k}
\ge 0$.
\end{lemma}
To prove Theorems~\ref{theorem_AD_particular}
and~\ref{theorem_Leger_particular} we first obtain sharp pointwise
estimates for the permutations related to the
kernels~(\ref{kappa_n}).
\begin{lemma}
\label{lemma_main_inequality} It holds that
\begin{equation}
\label{main_inequality}
 p_{\kappa_2}(z_1,z_2,z_3)\le 2 p_{\kappa_1}(z_1,z_2,z_3)\qquad \text{for all} \qquad
 (z_1,z_2,z_3)\in \mathbb{C}^3.
\end{equation}
\end{lemma}
\begin{Proof}
It is enough to prove (\ref{main_inequality}) for
$(z_1,z_2,z_3)=(0,u,v)$ as the permutations of the form
(\ref{permutation_Main}) are invariant under translations. Given
$u=(x,y)$ and $v=(a,b)$, by (\ref{representation_CMPT}) we get
\begin{align*}
2p_{\kappa_1}(0,&u,v)-p_{\kappa_2}(0,u,v)\\
&=\frac{2h_1(u,v)}{|u|^2|v|^2|u-v|^2}-\frac{2\,x^2a^2(x-a)^2\,h_1(u,v)+h_2(u,v)}{|u|^{4}|v|^{4}|u-v|^{4}}\\
&=\frac{2\left[|u|^2|v|^2|u-v|^2-\,x^2a^2(x-a)^2\right]\,h_1(u,v)-h_2(u,v)}{|u|^{4}|v|^{4}|u-v|^{4}}.
\end{align*}
Now we obtain a lower estimate of the expression in the square brackets
before $h_1(u,v)$. Expanding $|u|^2|v|^2|u-v|^2$ gives
\begin{align*}
&(x^2+y^2)(a^2+b^2)\left((x-a)^2+(y-b)^2\right)-x^2a^2(x-a)^2\\
&=x^2a^2(y-b)^2+(x^2b^2+a^2y^2+b^2y^2)\left((x-a)^2+(y-b)^2\right)\\
&\ge x^2a^2(y-b)^2+(x^2b^2+a^2y^2)(x-a)^2.
\end{align*}
Thus,
\begin{align*}
2p_{\kappa_1}(0,&u,v)-p_{\kappa_2}(0,u,v)\ge\frac{G(x,y,a,b)}{|u|^{4}|v|^{4}|u-v|^{4}},
\end{align*}
where
 $$
G(x,y,a,b):=
2(x^2a^2(y-b)^2+(x^2b^2+a^2y^2)(x-a)^2)\,h_1(u,v)-h_2(u,v).
 $$
 Notice that by Lemma~\ref{lemma_CMPT},
\begin{align*}
 h_1(u,v)&=ax(y-b)^{2}+x(x-a)b^{2}+a(a-x)y^{2}\\
 h_2(u,v)&=(ax)^{3}(y-b)^{4}+(x(x-a))^{3}b^{4}+(a(a-x))^{3}y^{4}.
\end{align*}
Consequently, to prove the required inequality it is enough to show
that $G(x,y,a,b)\ge 0$. We separate the discussion into three cases.

1) Let $a=0$. Then
$$
G(x,y,0,b)=2x^4b^2\cdot x^2b^2-x^6b^4=x^6b^4\ge 0.
$$

2) Let $b=0$. Then
\begin{align*}
G&(x,y,a,0)\\
&=2\left(a^2x^2y^2+a^2y^2(x-a)^2\right)\left(axy^2+a(a-x)y^2\right)-\left(a^3x^3y^4+a^3(a-x)^3y^4\right)\\
&=2a^3y^4\left(x^2+(x-a)^2\right)\left(x+(a-x)\right)-a^3y^4\left(x^3-(x-a)^3\right)\\
&=a^4y^4\left(2(x^2+(x-a)^2)-(x^2+x(x-a)+(x-a)^2)\right)\\
&=a^4y^4\left(x^2 -x(x-a)+(x-a)^2\right)\\
&=a^4y^4\left(x^2-ax+a^2\right)\\
&=a^4y^4\left((x-\tfrac{1}{2}a)^2+\tfrac{3}{4}a^2\right)\ge 0.
\end{align*}

3) Let $a\neq 0$ and $b\neq 0$. We divide $G(x,y,a,b)$ by
$a^6b^4$, put $\alpha=x/a$ and $\beta=y/b$ and take into account
that by Lemma~\ref{lemma_CMPT} in these settings one has
$$
\frac{h_k(u,v)}{a^{4k-2}b^{2k}}=\alpha^{2k-1}(\beta-1)^{2k}+\alpha^{2k-1}(\alpha-1)^{2k-1}-(\alpha-1)^{2k-1}\beta^{2k},\qquad
k=1,2.
$$
Therefore
\begin{align*}
&\frac{G(x,y,a,b)}{a^6b^4}\\
&\quad =2\left(\alpha^2(\beta-1)^2+(\alpha^2+\beta^2)(\alpha-1)^2\right)\left(\alpha(\beta-1)^2+\alpha(\alpha-1)-(\alpha-1)\beta^2\right)\\
&\quad\quad\quad-\left(\alpha^3(\beta-1)^4+\alpha^3(\alpha-1)^3-(\alpha-1)^3\beta^4\right).
\end{align*}
Removing brackets and further collecting terms give
\begin{align*}
\frac{G(x,y,a,b)}{a^6b^4}
&=\left(\alpha^2-\alpha+1\right)\left(\beta^4-4\alpha\beta^3+6\alpha^2\beta^2-4\alpha^3\beta+\alpha^4\right)\\
&=\left((\alpha-\tfrac{1}{2})^2+\tfrac{3}{4}\right)(\alpha-\beta)^4\ge
0.
\end{align*}

Thus  $G(x,y,a,b)$ is non-negative in all the cases and so we
are done.
\end{Proof}
\begin{remark}
\label{remark2}
 The inequality (\ref{main_inequality}) is sharp as
it is known from \cite[Lemma 3]{Chunaev2016} that
$$
2\left[\tfrac{\Re (z_1-z_2)}{|z_1-z_2|}\tfrac{\Re
(z_1-z_3)}{|z_1-z_3|}\tfrac{\Re
(z_2-z_3)}{|z_2-z_3|}\right]^{2}p_{\kappa_1}(z_1,z_2,z_3)\le
p_{\kappa_2}(z_1,z_2,z_3).
$$
Indeed, when all sides of the triangle $(z_1,z_2,z_3)$ make a small angle with the horizontal, the multiplier in the square brackets is close to $1$ in modulus.
\end{remark}
The estimate (\ref{main_inequality}) allows us to obtain an
inequality for $L^2$-norms.
\begin{lemma}
\label{lemma_inequality_for_norms} Let $\mu$ have linear growth and
$\varepsilon>0$ be fixed. Then
\begin{equation}
\label{L^2-norms} \|T_{\kappa_2,\varepsilon}1\|_{L^2(\mu)}\le
\sqrt{2}\,\|T_{\kappa_1,\varepsilon}1\|_{L^2(\mu)}+C\sqrt{\mu(\complex)},\qquad
C>0.
\end{equation}
\end{lemma}
\begin{Proof}
From Lemma~\ref{lemma_main_inequality} and the definition
(\ref{p_K(mu)}) of $p_K(\mu)$ we immediately get that
\begin{equation}
\label{inequality_perm_mu} p_{\kappa_2,\varepsilon}(\mu)\le 2
p_{\kappa_1,\varepsilon}(\mu).
\end{equation}
Now we use the identity
\begin{equation}
\label{Melnikov-Verdera-general}
\|T_{K,\varepsilon}1\|_{L^2(\mu)}^2=\tfrac{1}{3}p_{K,\varepsilon}(\mu)+\mathcal{R}_{K,\varepsilon}(\mu),\qquad
|\mathcal{R}_{K,\varepsilon}(\mu)|\le C_K\mu(\complex),\qquad C_K>0,
\end{equation}
where $K$ is any real antisymmetric Calder\'{o}n-Zygmund kernel with non-negative permutations.
This identity is a~generalization of (\ref{Melnikov-Verdera}) and is
  contained in \cite[Lemma 3.3]{CMPT} (see also \cite[Section
5]{Chunaev2016}). In these terms the inequality (\ref{inequality_perm_mu})
gives
$$
\tfrac{1}{3}p_{\kappa_2,\varepsilon}(\mu)+\mathcal{R}_{\kappa_2,\varepsilon}(\mu)\le
2
\left(\tfrac{1}{3}p_{\kappa_1,\varepsilon}(\mu)+\mathcal{R}_{\kappa_1,\varepsilon}(\mu)\right)+\mathcal{R}_{\kappa_2,\varepsilon}(\mu)-2\mathcal{R}_{\kappa_1,\varepsilon}(\mu),
$$
and, consequently,
$$
\|T_{\kappa_2,\varepsilon}1\|_{L^2(\mu)}^2\le
2\|T_{\kappa_1,\varepsilon}1\|_{L^2(\mu)}^2+C\mu(\complex),\qquad C>0.
$$
Applying the inequality $\sqrt{ax^2+b}\le\sqrt{a}x+\sqrt{b}$ valid
for $a,b,x\ge 0$, we get (\ref{L^2-norms}).
\end{Proof}
\begin{remark}
Lemmas~\ref{lemma_main_inequality}
and~\ref{lemma_inequality_for_norms} are particular cases of
\cite[Lemma 7]{CMPTcap} and \cite[Main Lemma]{Tolsa2004},
correspondingly, but with an explicit constant. The explicitness of
the constant is very important here and actually enables us to
obtain the result.
\end{remark}

We are now ready to prove Theorems~\ref{theorem_AD_particular}
and~\ref{theorem_Leger_particular}.

By (\ref{L^2-norms}) and the triangle inequality,
\begin{align*}
\|T_{k_t,\varepsilon}1\|_{L^2(\mu)}&=
\|(T_{\kappa_2,\varepsilon}+t\cdot
T_{\kappa_1,\varepsilon})1\|_{L^2(\mu)}\\
&\ge|t|\|T_{\kappa_1,\varepsilon}1\|_{L^2(\mu)}-\|T_{\kappa_2,\varepsilon}1\|_{L^2(\mu)}\\
&\ge
(|t|-\sqrt{2})\|T_{\kappa_1,\varepsilon}1\|_{L^2(\mu)}-C\sqrt{\mu(\complex)}.
\end{align*}
Consequently,
\begin{equation}
\label{norms} \|T_{\kappa_1,\varepsilon}1\|_{L^2(\mu)}\le
\frac{\|T_{k_t,\varepsilon}1\|_{L^2(\mu)}+C\sqrt{\mu(\complex)}}{|t|-\sqrt{2}},\qquad
|t|>\sqrt{2},
\end{equation}
and therefore for any cube $Q \subset \mathbb{C}$,
$$
\|T_{\kappa_1,\varepsilon}\chi_Q\|_{L^2(\mu\lfloor Q)}\le
\frac{\|T_{k_t,\varepsilon}\chi_Q\|_{L^2(\mu\lfloor
Q)}+C\sqrt{\mu(Q)}}{|t|-\sqrt{2}},\qquad |t|>\sqrt{2}.
$$
Applying a variant of the $T1$ Theorem of Nazarov, Treil and Volberg
from \cite[Theorem 9.42]{Tolsa_book}, we infer that the
$L^2(\mu)$-boundedness of $T_{k_t}$ with $|t|>\sqrt{2}$ implies that
$T_{\kappa_1}$ (and hence the Cauchy transform) is
$L^2(\mu)$-bounded. Therefore, by Theorems~\ref{theorem_AD_MMV}
and~\ref{theorem_Leger}, we get the desired result. Note that the
``only if'' part of Theorem~\ref{theorem_AD_particular} follows from
\cite{David}.

\begin{remark}
\label{conjecture-inequality}
Computer experiments suggest that the following inequality holds:
\begin{equation}
\label{conj_ineq} p_{\kappa_N}(z_1,z_2,z_3)\le \tfrac{N}{n}
p_{\kappa_n}(z_1,z_2,z_3).
\end{equation}
(Lemma~\ref{lemma_main_inequality} corresponds to the case
$(n,N)=(1,2)$.) Moreover, if $u=-\gamma+i$, $v=\gamma+i$ and
$\gamma>0$, then (see \cite[Example 1]{Chunaev2016})
$$
p_{\kappa_m}(0,u,v)=\frac{\gamma^{2m-2}\left((\gamma^2+1)^m-\gamma^{2m}\right)}{(\gamma^2+1)^{2m}},\qquad
m\in \mathbb{N},
$$
and therefore
$$
\lim_{\gamma\to\infty}\frac{p_{\kappa_N}(0,u,v)}{p_{\kappa_n}(0,u,v)}=
\lim_{\gamma\to\infty}\frac{1-\left(\gamma^2/(\gamma^2+1)\right)^N}{1-\left(\gamma^2/(\gamma^2+1)\right)^n}
=\frac{N}{n}.
$$
It means that the constant $N/n$ is sharp if (\ref{conj_ineq}) is
true.

It would follow from (\ref{conj_ineq}) in the same manner as above that
the $L^2(\mu)$-boundedness of $T_{K_t}$ with $|t|>\sqrt{N/n}$ implies that $T_{\kappa_n}$ is $L^2(\mu)$-bounded. This would give the analogues of Theorems~\ref{theorem_AD_particular} and~\ref{theorem_Leger_particular}
for the more general case of kernels~(\ref{kernels_n_N}) via theorems in \cite{CMPT} instead of Theorems~\ref{theorem_AD_MMV} and~\ref{theorem_Leger}.
\end{remark}

\section{The proof of Theorem~\ref{theorem_AD_general}}
\label{section_AD_general}

We now come to the case of kernels~(\ref{kernels_n_N}) with $t\in
\Omega(n,N)$ (see (\ref{Omega})). As mentioned above, the
permutations are non-negative for them and hence a curvature-like
method can be used directly. Namely, we will adapt the arguments
from \cite[Section~8]{CMPT} which in turn stem from
\cite{DavSem_Ast} to our settings. Note that the ``only if'' part of
Theorem~\ref{theorem_AD_general} follows from \cite{David} (even for
all $t\in \mathbb{R}$). Thus we only need to prove the ``if'' part.

Suppose that $\mu$ is an AD-regular measure on $\mathbb{C}$ and
$T_{K_t}$ the CZO associated with the kernels $(\ref{kernels_n_N})$,
$t\in \Omega(n,N)$. It is proved in \cite[Lemmas 5 and
6]{Chunaev2016} that if
\begin{eqnarray}
\label{t1}&t\in \mathbb{R}\setminus \left[-\tfrac{1}{2}\left(3+\sqrt{9-4\tfrac{N}{n}}\right); 2-\tfrac{N}{n}\right],\qquad & n<N\le 2n,\\
\label{t2}&t\in \mathbb{R}\setminus
\left[-\frac{1}{2}\left(\sigma_{n,M}+\sqrt{\sigma_{n,M}^2-4\frac{N}{n}}\right);
\sigma_{n,M}-3\right],\qquad &N\ge 2n,
\end{eqnarray}
where $\sigma_{n,M}=3+\left(\tfrac{N}{n}-2\right)\sqrt{N-2n}$ as above, then
$$
p_{K_t}(z_1,z_2,z_3)\ge C(t)\cdot p_{\kappa_n}(z_1,z_2,z_3),\qquad C(t)>0,\qquad (z_1,z_2,z_3)\in \mathbb{C}^3.
$$
Consequently, $p_{K_t,\varepsilon}(\mu)\ge C(t)\cdot
p_{\kappa_n,\varepsilon}(\mu)$  and
 hence from (\ref{Melnikov-Verdera-general}) we
 conclude that for $t$ as in (\ref{t1}) and (\ref{t2}) and any cube $Q\subset \mathbb{C}$,
\begin{equation}
\label{norms_K_t}
\|T_{\kappa_n,\varepsilon}\chi_Q\|_{L^2(\mu\lfloor Q)}\le
C(t)\left(\|T_{K_t,\varepsilon}\chi_Q\|_{L^2(\mu\lfloor Q)}+C\sqrt{\mu(Q)}\right).
\end{equation}
Therefore, by the $T1$ Theorem from \cite[Theorem 9.42]{Tolsa_book} and \cite[Theorem
1.3]{CMPT}, the measure $\mu$ is uniformly rectifiable.

What is left, according to (\ref{Omega}), is to prove Theorem~\ref{theorem_AD_general} for
\begin{eqnarray}
& \label{end_1} t=2-\tfrac{N}{n},\qquad & n<N\le 2n,\phantom{\sqrt{\tfrac{N}{n}}}\\
&\label{end_2}
t=\sigma_{n,M}-3,\phantom{\sqrt{\tfrac{N}{n}}}&N\ge 2n,\\
& \label{end_3} t=-\tfrac{1}{2}\left(3+\sqrt{9-4\tfrac{N}{n}}\right), & n<N\le 2n,\\
&\label{end_4}
t=-\frac{1}{2}\left(\sigma_{n,M}+\sqrt{\sigma_{n,M}^2-4\frac{N}{n}}\right),\qquad&N\ge
2n.
\end{eqnarray}

To manage these cases, we introduce additional notation. Given two
distinct points $z,w\in \mathbb{C}$, we denote by $L_{z,w}$ the line
passing through $z$ and $w$. Given three pairwise distinct points
$z_1,z_2,z_3\in \mathbb{C}$, we denote by
$\measuredangle(z_1,z_2,z_3)$ the smallest angle (belonging to
$[0;\pi/2]$) formed by the lines $L_{z_1,z_2}$ and $L_{z_1,z_3}$. If
$L$ and $L'$ are lines, then $\measuredangle(L,L')$ is the smallest
angle (belonging to $[0;\pi/2]$) between them. Also,
$\theta_V(L):=\measuredangle(L,V)$, where $V$ is the vertical.
Furthermore, for a fixed constant $\tau\ge 1$ and complex numbers
$z_1$, $z_2$ and $z_2$, set
\begin{equation}
\label{O_tau}
\mathcal{O}_\tau:=\left\{(z_1,z_2,z_3):\frac{|z_i-z_j|}{|z_i-z_k|}\le\tau
\text{ for pairwise distinct } i,j,k\in \{1,2,3\}\right\},
\end{equation}
so that all triangles with vertexes $z_1$, $z_2$ and $z_3$ in $\mathcal{O}_\tau$ have comparable sides.

Given $\alpha_0\in(0,\pi/2)$ and $(z_1,z_2,z_3)$, in what follows sometimes we will need the conditions
\begin{equation}
\label{delta1}
\theta_{V}(L_{z_1,z_2})+\theta_{V}(L_{z_2,z_3})+\theta_{V}(L_{z_1,z_3})\ge
\alpha_0
\end{equation}
and
\begin{equation}
\label{delta2}
\theta_{V}(L_{z_1,z_2})+\theta_{V}(L_{z_2,z_3})+\theta_{V}(L_{z_1,z_3})\le
\tfrac{3}{2}\pi-\alpha_0.
\end{equation}
We will also use the following result.
\begin{lemma}[Lemma 10 in \cite{Chunaev2016}]
\label{lemma_ends} Fix $\alpha_0\in(0,\pi/2)$. Given $K_t$ and
$(z_1,z_2,z_3)\in \mathcal{O}_\tau$,

 $($i$)$ if the condition
$(\ref{delta1})$ is satisfied and $t$ is as in $(\ref{end_1})$ or
  $(\ref{end_2})$,

or

$($ii$)$ if both the conditions $(\ref{delta1})$ and
$(\ref{delta2})$ are satisfied and $t$ is as in $(\ref{end_3})$ or
  $(\ref{end_4})$,

then the following inequality holds
\begin{equation}
\label{last_section_ineq}
p_{K_t}(z_1,z_2,z_3)\ge C(\alpha_0,\tau)\cdot c(z_1,z_2,z_3)^2,\qquad C(\alpha_0,\tau)>0.
\end{equation}
\end{lemma}
On the one hand, if we are in the clause $(i)$ of
Lemma~\ref{lemma_ends}, i.e. in the same settings as in~\cite{CMPT}, then we
can undeviatingly follow the scheme from \cite[Section 8]{CMPT}
(exchanging $p_{\kappa_n}$ for $p_{K_t}$) in order to get our result
for $t$ as in~(\ref{end_1}) or~(\ref{end_2}).

On the other hand, by the clause $(ii)$ of Lemma~\ref{lemma_ends}, we can ensure that
the inequality~(\ref{last_section_ineq}) is true for $t$ as in
$(\ref{end_3})$ or $(\ref{end_4})$ if the sides of the
triangles $(z_1,z_2,z_3)$ are far from both the vertical and
horizontal. Consequently,  the scheme from \cite[Section
8]{CMPT} cannot be applied directly for such $t$. Nevertheless, as
we show below, it works after several modifications (besides the
exchange of $p_{\kappa_n}$ for $p_{K_t}$) connected basically with
adapting geometrical arguments to both the conditions
$(\ref{delta1})$ and $(\ref{delta2})$. Note that some of the arguments in~\cite[Section 8]{CMPT} are very sketchy and so, for the sake of completeness, we give a proof that is more detailed than the corresponding one in \cite[Section 8]{CMPT}.

The fact that the $L^2(\mu)$-boundedness of $T_{K_t}$ implies that
$\mu$ is uniformly rectifiable will be proved by means of a corona
type decomposition. We now recall how such a decomposition is
defined in \cite[Chapter 2]{DavSem_Ast} for a given AD-regular Borel
measure $\mu$. The elements $Q$ playing the role of dyadic cubes are
usually called $\mu$-cubes.

Given a $1$-dimensional AD regular Borel measure $\mu$ on
$\mathbb{C}$, for each $j\in \mathbb{Z}$ (or $j\ge j_0$ if
$\mu(\mathbb{C})<\infty$) there exists a family $\mathcal{D}_j$ of
Borel subsets of $\textrm{spt}\,\mu$, i.e. $\mu$-cubes $Q$ of the
$j$th generation, such that:
\begin{enumerate}
  \item each $\mathcal{D}_j$ is a disjoint partition of $\textrm{spt}\;\mu$,
  i.e. if $Q,Q'\in \mathcal{D}_j$ and $Q\neq Q'$, then
  $$
  \textrm{spt}\,\mu=\bigcup_{Q\in \mathcal{D}_j}Q\quad \text{and}\quad Q\cap
  Q'=\varnothing;
  $$
  \item if $Q\in \mathcal{D}_j$ and $Q'\in \mathcal{D}_k$ with $k\le j$,
   then either $Q\subseteq  Q'$ or $Q\cap Q'=\varnothing$;
  \item for all $j\in \mathbb{Z}$ and $Q\in D_j$, we have
  $$
  2^{-j}\lesssim \diam(Q)\lesssim 2^{-j}\quad \text{and} \quad \mu(Q)\approx
  2^{-j}.
  $$
\end{enumerate}

In what follows, $\mathcal{D}:=\bigcup_{j\in
\mathbb{Z}}\mathcal{D}_j$. Moreover, given $Q\in \mathcal{D}_j$, we define the
\textit{side length} of $Q$ as $\ell(Q)=2^{-j}$, which actually
indicates the generation of $Q$. Obviously,
$\ell(Q)\approx\diam(Q)$. The value of $\ell(Q)$ is not well defined
if the $\mu$-cube $Q$ belongs to $\mathcal{D}_j\cap \mathcal{D}_k$
with $j\neq k$. To avoid this, one may consider a $Q\in \mathcal{D}_j$ as a couple $(Q,j)$.

Given $\lambda>1$ and $Q\in \mathcal{D}$, set
$$
\lambda Q:=\{x\in \textrm{spt}\,\mu: \dist(x,Q)\le (\lambda-1)\ell(Q)\}.
$$
We will also need the following version of P. Jones' $\beta$-numbers
for $\mu$-cubes (see \cite{David_Semmes}):
\begin{equation*}
\label{beta-cubes}
\beta_q(Q)=\inf_{L}\left(\frac{1}{\ell(Q)}\int_{\eta_1 Q}\left(\frac{\dist(x,L)}{\ell(Q)}\right)^q
d\mu(x)\right)^{1/q},\qquad 1\le q\le \infty,
\end{equation*}
where $\eta_1>4$ is some fixed constant and the infimum is taken over all affine lines $L$. We will mostly
use $\beta_1(Q)$ and denote by $L_Q$ the best approximating
line for $\beta_1(Q)$.

Given $Q\in \mathcal{D}_j$, the \textit{sons} of $Q$, forming the collection $\textsf{Sons}(Q)$, are the $\mu$-cubes
$Q'\in \mathcal{D}_{j+1}$ such that $Q'\subseteq Q$.

By \cite[Chapter 2]{DavSem_Ast}, one says that $\mu$ admits a
\textit{corona decomposition} if there are $\eta>0$, $\delta>0$ and a triple $(\mathcal{B},\mathcal{G},{\sf Tree})$, where
$\mathcal{B}$ and $\mathcal{G}$ are subsets of $\mathcal{D}$ (the
``bad $\mu$-cubes'' and the ``good $\mu$-cubes'') and ${\sf Tree}$
is a family of subsets $S$ of $\mathcal{G}$ so that the
following conditions are satisfied:
\begin{enumerate}
  \item $\mathcal{D}=\mathcal{B}\cup \mathcal{G}$ and $\mathcal{B}\cap
  \mathcal{G}=\varnothing$.
  \item $\mathcal{B}$ satisfies a Carleson packing condition, i.e.
  \begin{equation}
  \label{Carleson}
  \sum_{Q\in \mathcal{B}: Q\subseteq R}\mu(Q)\lesssim_\eta \mu(R)
  \quad\text{for all} \quad R\in \mathcal{D}.
  \end{equation}
  \item $\mathcal{G}=\bigcup_{S\in {\sf Tree}}S$ and the union is disjoint;
  \item Each $S\in {\sf Tree}$ is called a \textit{tree} and is \textit{coherent}: each $S$ has
  a  unique maximal element $Q_S$, which contains all other elements
  of $S$ as subsets, i.e.
\begin{itemize}
  \item a $\mu$-cube $Q'\in \mathcal{D}$ belongs to $S$ if $Q\subseteq Q'\subseteq
  Q_S$ for some $Q\in S$;
  \item if $Q\in S$ then either all elements of $\textsf{Sons}(Q)$ lie in
  $S$ or none of them do.
\end{itemize}
  \item For each $S\in {\sf Tree}$, there exists a (possibly
  rotated) Lipschitz graph $\Gamma_S$ with constant smaller than
  $\eta$ such that $\dist(x,\Gamma_S)\le \delta \diam(Q)$ whenever $x\in
  2Q$ and $Q\in S$.
  \item The maximal $\mu$-cubes $Q_S$, for $S\in {\sf Tree}$,
  satisfy the Carleson packing condition
  $$
\sum_{S\in {\sf Tree}:\; Q_S\subseteq R}\mu(Q_S)\lesssim \mu(R)
  \quad\text{for all} \quad R\in \mathcal{D}.
  $$
\end{enumerate}

According to \cite{DavSem_Ast} (see e.g. Section~1, (C4) and (C6)),
if $\mu$ is uniformly rectifiable, then it admits a corona
decomposition for all parameters $\eta,\delta>0$. Conversely, the
existence of a corona decomposition for a single set of the
parameters $\eta$ and $\delta$ implies that $\mu$ is uniformly
rectifiable.

We now turn to constructing a corona decomposition for our settings.
From now on, $\textsf{B}(\varepsilon)$ stands for the family of
cubes $Q\in \mathcal{D}$ such that $\beta_1(Q)\ge \varepsilon$.
Furthermore, $\textsf{G}(\varepsilon):=\mathcal{D}\setminus
\textsf{B}(\varepsilon)$. The aim  is to show that
$\textsf{B}(\varepsilon)$ satisfies a Carleson packing condition.

By H\"{o}lder's inequality, $\beta_{1}(Q)\lesssim \beta_{2}(Q)$. Thus, for any ${\varepsilon>0}$, if
$\beta_{1}(Q)\ge \varepsilon$, then $\beta_{2}(Q)\gtrsim
\varepsilon$. Moreover, for any $y,z \in \eta_1 Q$, we have
\begin{align*}
\beta_{2}(Q)^2&\le
\frac{1}{\ell(Q)}\int_{\eta_1Q}\bigg(\frac{\dist(x,L_{y,z})}{\ell(Q)}\bigg)^2d\mu(x)\\
& =\frac{1}{\ell(Q)}\bigg(\int_{\tiny \begin{array}{l}
                                x\in\eta_1Q: \\
                                \frac{\dist(x,L_{y,z})}{\ell(Q)}<\varepsilon^2
                              \end{array}}\bigg(\frac{\dist(x,L_{y,z})}{\ell(Q)}\bigg)^2d\mu(x)\\
                              &\qquad\qquad\qquad +
                              \int_{\tiny \begin{array}{l}
                                x\in\eta_1Q: \\
                                \frac{\dist(x,L_{y,z})}{\ell(Q)}\ge \varepsilon^2
                              \end{array}}\bigg(\frac{\dist(x,L_{y,z})}{\ell(Q)}\bigg)^2d\mu(x)\bigg)\\
                              &\lesssim
                              \frac{1}{\ell(Q)}\bigg(\varepsilon^4\ell(Q)+
                              \int_{\tiny \begin{array}{l}
                                x\in\eta_1Q: \\
                                \frac{\dist(x,L_{y,z})}{\ell(Q)}\ge \varepsilon^2
                              \end{array}}\bigg(\frac{\dist(x,L_{y,z})}{\ell(Q)}\bigg)^2d\mu(x)\bigg)\\
                              &=
                              \varepsilon^4 +
                              \frac{1}{\ell(Q)}\int_{\tiny \begin{array}{l}
                                x\in\eta_1Q: \\
                                \frac{\dist(x,L_{y,z})}{\ell(Q)}\ge \varepsilon^2
                              \end{array}}\bigg(\frac{\dist(x,L_{y,z})}{\ell(Q)}\bigg)^2d\mu(x).
\end{align*}
We used that $\ell(Q)\approx \mu(Q)$.

\begin{lemma}
\label{lemma_new}
Let $B_1=B(\zeta_1,r_1)$ and $B_2=B(\zeta_2,r_2)$ be two balls such that $B_1\subset \eta_1Q$, $B_2\subset \eta_1Q$, $\dist (B_1,B_2)\approx \ell(Q)$, $\zeta_1,\zeta_2\in \eta_1Q$ and $r_1\approx r_2\approx \ell(Q)$. If $y\in B_1$ and $z\in B_2$, then for $\varepsilon$ small enough,
$$\int_{\tiny \begin{array}{l}
                                x\in\eta_1Q: \\
                                \frac{\dist(x,L_{y,z})}{\ell(Q)}\ge \varepsilon^2
                              \end{array}}\dist(x,L_{y,z})^2d\mu(x)\lesssim_\varepsilon
                               \ell(Q)^2p_{K_t}^{(\varepsilon;\,Q)}(\mu),
$$
where
\begin{equation*}
\label{p_Q(mu)}
 p_{K_t}^{(\varepsilon;\,Q)}(\mu):=\iiint _{\tiny
\begin{array}{l}
          (x,y,z)\in (\eta_1Q)^3 \\
          |x-y|\ge\varepsilon^2 \ell(Q),\\
          |x-z|\ge\varepsilon^2 \ell(Q)
        \end{array}
} p_{K_t}(x,y,z)\;d\mu(x)d\mu(y)d\mu(z).
\end{equation*}
\end{lemma}

Note that the existence of the above-mentioned balls $B_1$ and $B_2$ is guaranteed in the AD-regular case.

\smallskip

\begin{Proof}
First note that $|x-y|\ge\varepsilon^2 \ell(Q)$ and $|x-z|\ge\varepsilon^2 \ell(Q)$ as $\dist(x,L_{y,z}) \ge \varepsilon^2\ell(Q)$. Consequently, since $x\in \eta_1 Q$,
$y\in B_1$ and $z\in B_2$,
\begin{equation}
\label{comparability}
|x-z|\approx |x-y|\approx |y-z|,
\end{equation}
where the comparability constants depend on $\eta_1$ and $\varepsilon$.
We now separate two cases.

$(1)$ Suppose that
$$
\varepsilon^{10}\le \theta_{V}(L_{x,y})+\theta_{V}(L_{y,z})+\theta_{V}(L_{x,z})\le
\tfrac{3}{2}\pi-\varepsilon^{10}.
$$
Then by the clause $(ii)$ of Lemma~\ref{lemma_ends}, where we put $\alpha_0=\varepsilon^{10}$ and $\tau=\tau(\varepsilon,\eta_1)$ chosen from (\ref{comparability}), we have $c(x,y,z)^2\lesssim_\varepsilon p_{K_t}(x,y,z)$.

$(2)$ Now let
$$
\theta_{V}(L_{x,y})+\theta_{V}(L_{y,z})+\theta_{V}(L_{x,z})<\varepsilon^{10}
$$
or
$$
\theta_{V}(L_{x,y})+\theta_{V}(L_{y,z})+\theta_{V}(L_{x,z})>
\tfrac{3}{2}\pi-\varepsilon^{10}.
$$
In this case $\dist (x, L_{y,z})\lesssim \varepsilon^{10} \ell(Q)$. Thus for $\varepsilon$ small enough we get a contradiction with the assumption $\dist(x,L_{y,z}) \ge \varepsilon^2 \ell(Q)$.

Summarizing,
\begin{align*}
&\int_{\tiny \begin{array}{l}
                                x\in\eta_1Q: \\
                                \frac{\dist(x,L_{y,z})}{\ell(Q)}\ge \varepsilon^2
                              \end{array}}\dist(x,L_{y,z})^2d\mu(x)\\
&\qquad \lesssim \frac{\ell(Q)^4}{\mu(B_1)\mu(B_2)}
 \int_{B_2}\int_{B_1}\int_{\tiny \begin{array}{l}
                                \eta_1Q: \\
                                \frac{\dist(x,L_{y,z})}{\ell(Q)}\ge \varepsilon^2
                              \end{array}}\left(\frac{\dist(x,L_{y,z})}{|x-y||x-z|}\right)^2d\mu(x) d\mu(y)d\mu(z)\\
&\qquad \lesssim  \ell(Q)^2
 \int_{\eta_1Q}\int_{\eta_1Q}\int_{\tiny \begin{array}{l}
                                \eta_1Q: \\
                                \frac{\dist(x,L_{y,z})}{\ell(Q)}\ge \varepsilon^2
                              \end{array}}c(x,y,z)^2d\mu(x) d\mu(y)d\mu(z)\\
&\qquad \lesssim_{\varepsilon}  \ell(Q)^2
 p_{K_t}^{(\varepsilon;\,Q)}(\mu).
\end{align*}
We used the well-known identity
$c(x,y,z)=\dist(x,L_{y,z})/(|x-y||x-z|)$.
\end{Proof}

The estimate for $\beta_2(Q)^2$ that we obtained above and Lemma~\ref{lemma_new} give
$$
\beta_2(Q)^2\lesssim \varepsilon^4+\frac{C(\varepsilon)}{\ell(Q)}\;
 p_{K_t}^{(\varepsilon;\,Q)}(\mu),
$$
so, taking into account that $\beta_{2}(Q)\gtrsim
\varepsilon$, we get for sufficiently small $\varepsilon$  that
\begin{equation*}
\label{mu(Q)}
 \mu(Q)\lesssim_{\varepsilon}p_{K_t}^{(\varepsilon;\,Q)}(\mu) \qquad \text{for all} \qquad Q\in \mathcal{D}\quad \text{such that} \quad \beta_{1}(Q)\ge
 \varepsilon.
\end{equation*}
From this we immediately get that
$$
\sum_{Q\in \textsf{B}(\varepsilon):\, Q\subseteq R}\mu(Q)
\lesssim_\varepsilon \sum_{Q\in \textsf{B}(\varepsilon):\, Q\subseteq
R} p_{K_t}^{(\varepsilon;Q)}(\mu).
$$
To estimate the latter sum, we will use the notation
\begin{equation*}
\label{annuli} A_j(\varepsilon):=\{x:\varepsilon^2\ell(Q) \le
|x-y|\le C\ell(Q)\}, \qquad Q\in \mathcal{D}_j,\qquad C>0.
\end{equation*}
These are the concentric annuli $B(y,C\ell(Q))\setminus
B(y,\varepsilon^2 \ell(Q))$, contained in the ball $B(y,C\ell(R))$,
where $\ell(R)=2^{-j_0}$. They have bounded overlap  depending on $\varepsilon$ and thus the sum
$\sum_{j\ge j_0} \int_{R\cap A_j(\varepsilon)}$ is less than
$C(\varepsilon)\int_{R}$ with some $C(\varepsilon)>0$. These
observations lead to the following:
\begin{align*}
&\sum_{Q\in \textsf{B}(\varepsilon):\, Q\subseteq R}
p_{K_t}^{(\varepsilon;Q)}(\mu)\\
&\qquad =\sum_{j\ge j_0}\sum_{Q\in
\textsf{B}(\varepsilon)\cap \mathcal{D}_j(R)}  \int_{\eta_1Q} \int_{\eta_1Q} \int_{\eta_1Q\cap A_j(\varepsilon)} p_{K_t}(x,y,z)\;d\mu(x)d\mu(y)d\mu(z)\\
&\qquad \le \int_{\eta_1R}\int_{\eta_1R} \Bigg(\sum_{j\ge j_0}\sum_{Q\in
\textsf{B}(\varepsilon)\cap \mathcal{D}_j(R)} \int_{\eta_1Q\cap A_j(\varepsilon)} p_{K_t}(x,y,z)\;d\mu(x)\Bigg)d\mu(y)d\mu(z)\\
&\qquad = \int_{\eta_1R}\int_{\eta_1R} \Bigg(\sum_{j\ge j_0} \int_{\eta_1R\cap A_j(\varepsilon)} p_{K_t}(x,y,z)\;d\mu(x)\Bigg)d\mu(y)d\mu(z)\\
&\qquad \lesssim_\varepsilon
\iiint_{(\eta_1R)^3} p_{K_t}(x,y,z)\;d\mu(x)d\mu(y)d\mu(z)\\
&\qquad = p_{K_t}(\mu \lfloor (\eta_1R)).
\end{align*}
Since $T_{K_t}$ is $L^2(\mu)$-bounded, we get $p_{K_t}(\mu\lfloor F)\lesssim
\mu(F)$ for any $F\subset \mathbb{C}$. Consequently,
$$
p_{K_t}(\mu \lfloor (\eta_1R))\lesssim
\mu(R)\qquad \text{for all}\qquad R\in \mathcal{D},
$$
and therefore we reach the desired inequality
\begin{equation*}
\label{8.2} \sum_{Q\in \textsf{B}(\varepsilon):\, Q\subseteq
R}\mu(Q)\lesssim_\varepsilon\mu(R) \qquad \text{for all}\qquad R\in
\mathcal{D}.
\end{equation*}

Thus, for any $\varepsilon>0$, there exists the
decomposition
\begin{equation}
\label{initial_decomp}
 \mathcal{D}=\textsf{B}(\varepsilon)\cup
\textsf{G}(\varepsilon),
\end{equation}
where $\textsf{B}(\varepsilon)$ satisfies a Carleson packing
condition and for any cube $Q\in \textsf{G}(\varepsilon)$ there
exists a line $L_Q$ such that $\dist (x,L_Q)\lesssim \sqrt{\varepsilon}\,\ell(Q)$
for all $x\in \tfrac{1}{2}\eta_1Q$ (since $\beta_1(Q)<\varepsilon$ for such cubes and $\beta_\infty(Q)\lesssim \sqrt{\beta_1(2Q)}$). More details can be found
in \cite[Ch. 6]{DavSem_Ast}.

Using the decomposition (\ref{initial_decomp}), we now can apply
\cite[Lemma~7.1]{DavSem_Ast} in order to obtain a new decomposition
but already with a family of stopping times regions. Let
$\theta_0:=10^{-6}\pi$. Note that one has to choose $\varepsilon\ll
\theta_0$ to prove the following assertion.
\begin{lemma}
\label{lemma8.1} There exists a decomposition
$\mathcal{D}=\mathcal{B}\cup \mathcal{G}$, where $\mathcal{B}$
satisfies a Carleson packing condition and $\mathcal{G}$ can be
partitioned into a family {\rm Tree} of coherent regions $S$,
satisfying the following. For each $S\in {\rm Tree}$ denote
$$
\alpha(S):=\tfrac{1}{10}\theta_0 \quad \text{if} \quad
\theta_0\le\theta_V(L_{Q_S})\le\pi/2-\theta_0
$$
and
$$
\alpha(S):=10\theta_0 \quad \text{if} \quad \theta_V(L_{Q_S})<
\theta_0 \quad \text{or} \qquad \theta_V(L_{Q_S})>
\pi/2-\theta_0.
$$
Then we have
\begin{itemize}
  \item if $Q\in S$, then $\measuredangle(L_Q,L_{Q_S})\le \alpha(S)$;
  \item if $Q$ is a minimal cube of $S$, then either at least one element of ${\sf Sons}(Q)$ lies in $\mathcal{B}$ or else $\measuredangle(L_Q,L_{Q_S})\ge
  \tfrac{1}{2}\alpha(S)$.
\end{itemize}
\end{lemma}

Here $\mathcal{G}\subseteq \textsf{G}(\varepsilon)$ and
therefore for any $Q\in \mathcal{G}$ one has
$\beta_1(Q)<\varepsilon$.

Lemma~\ref{lemma8.1} is an analogue of \cite[Lemma 8.1]{CMPT} which
 comes from \cite[Lemma 7.1]{DavSem_Ast}. The main difference
between \cite[Lemma 8.1]{CMPT} and \cite[Lemma 7.1]{DavSem_Ast} is
that two different values of the parameter $\alpha(S)$ have to be
chosen, according to the angle $\theta_V(L_{Q_S})$. In our case  the
situations where the angle $\theta_V(L_{Q_S})$ is close to zero and
$\pi/2$ have to be also distinguished.

To obtain the required Lipschitz graph, one can follow the proof of
\cite[Proposition 8.2]{DavSem_Ast}. This leads to the following
statement.
\begin{lemma}
For each $S\in {\sf Tree}$ from Lemma~\ref{lemma8.1}, there exists a
Lipschitz function $A_S:L_{Q_S}\to L^\bot_{Q_S}$ with norm $\lesssim
 \alpha(S)$ such that, denoting by $\Gamma_S$ the graph of $A_S$,
$$
\dist(x, \Gamma_S)\lesssim \sqrt{\varepsilon}\, \ell(Q)
$$
for all $x\in 2Q$ with $Q\in S$.
\end{lemma}

The proof  will be completed if we show that the maximal $\mu$-cubes
$Q_S$, $S\in {\sf Tree}$, satisfy the Carleson packing condition
$$
\sum_{S\in \textsf{Tree}:\; Q_S\subseteq R}\mu(Q_S)\lesssim
\mu(R)\qquad \text{for all}\qquad R\in \mathcal{D}.
$$
To do so, we will distinguish several types of trees.

Here and subsequently, ${\sf Stop}(S)$ denotes the family of the
\textit{minimal} $\mu$-cubes of ${S\in {\sf Tree}}$, which may be
empty. By Lemma~\ref{lemma8.1}, we can split ${\sf Stop}(S)$ as
follows:
\begin{equation}
\label{Stop} {\sf Stop}(S)={\sf Stop}_\alpha(S) \cup {\sf
Stop}_\beta(S),\qquad {\sf Stop}_\alpha(S) \cap {\sf
Stop}_\beta(S)=\varnothing,
\end{equation}
where ${\sf Stop}_\beta(S)$ contains all minimal $\mu$-cubes $Q$
such that at least one element of $\textsf{Sons}(Q)$ belongs to $\mathcal{B}$,
and ${\sf Stop}_\alpha(S)$ contains all minimal  $Q$ such that
$\measuredangle(L_Q,L_{Q_S})\ge \tfrac{1}{2}\alpha(S)$.

The first set that we will consider is
$$
\Delta_1:=\left\{S\in \textsf{Tree}:\mu\left(Q_S\setminus
\medcup_{Q\in \textsf{Stop}(S)}Q\right)\ge
\tfrac{1}{2}\,\mu(Q_S)\right\}.
$$
Clearly, if $S\in \textsf{Tree}\setminus \Delta_1$, then by
(\ref{Stop}),
\begin{equation}
\label{not_I} \tfrac{1}{2}\,\mu(Q_S)< \mu\left(\medcup_{Q\in
\textsf{Stop}(S)}Q\right)=\mu\left(\medcup_{Q\in
\textsf{Stop}_\alpha(S)}Q\right)+\mu\left(\medcup_{Q\in
\textsf{Stop}_\beta(S)}Q\right).
\end{equation}
Now let
$$
\Delta_2:=\left\{S\in \textsf{Tree}\setminus
\Delta_1:\mu\left(\medcup_{Q\in \textsf{Stop}_\beta(S)}Q\right)\ge
\tfrac{1}{4}\,\mu(Q_S)\right\},
$$
i.e. not less than $\frac{1}{4}\,\mu(Q_S)$ of the measure of the
minimal cubes for these trees have sons in $\mathcal{B}$. The rest of the trees are in
$$
\Delta_3:=\left\{S\in \textsf{Tree}\setminus (\Delta_1\cup
\Delta_2):\mu\left(\medcup_{Q\in \textsf{Stop}_\alpha(S)}Q\right)\ge
\tfrac{1}{4}\,\mu(Q_S)\right\}.
$$
Indeed, if $S\in \textsf{Tree}\setminus (\Delta_1\cup \Delta_2\cup
\Delta_3)$, then (\ref{not_I}) is valid and moreover
$$
\mu\left(\medcup_{Q\in \textsf{Stop}_\alpha(S)}Q\right)<
\tfrac{1}{4}\,\mu(Q_S) \quad \text{and} \quad \mu\left(\medcup_{Q\in
\textsf{Stop}_\beta(S)}Q\right)< \tfrac{1}{4}\,\mu(Q_S).
$$
This means that $\textsf{Tree}\setminus (\Delta_1\cup \Delta_2\cup
\Delta_3)=\varnothing$.

We also split $\Delta_3$ in the three disjoint sets:
\begin{align*}
&\Delta^{'\phantom{''}}_{3} :=\left\{S\in
\Delta_3:\theta_0\le\theta_V(L_{Q_S})\le\pi/2-\theta_0\right\},\\
&\Delta^{''\phantom{'}}_{3}:=\left\{S\in
\Delta_3:\theta_V(L_{Q_S})< \theta_0\right\},\\
&\Delta^{'''}_{3}:=\left\{S\in
\Delta_3:\theta_V(L_{Q_S})> \pi/2-\theta_0\right\}.
\end{align*}

So we have the disjoint union
$$
\textsf{Tree}=\Delta_1\cup \Delta_2\cup \Delta^{'}_3\cup \Delta^{''}_3\cup \Delta^{'''}_3.
$$
The procedure now is to check the required Carleson packing
condition for all components of this union.

For all $S\in \textsf{Tree}$ the sets $Q_S\setminus \bigcup_{Q\in
\textsf{Stop}(S)}Q$ are pairwise disjoint and hence
$$
\sum_{S\in \Delta_1:\;Q_S\subseteq R}\mu(Q_S)\le 2\sum_{S\in
\textsf{Tree}:\;Q_S\subseteq R}\mu\left(Q_S\setminus \medcup_{Q\in
\textsf{Stop}(S)}Q\right)\le 2\mu(R).
$$
If $S\in \Delta_2$, then by definition and the fact that
$\mu(Q)\approx \mu(Q')$ for $Q'\in \textsf{Sons}(Q)$,
$$
\mu(Q_S)\le 4\mu\left(\medcup_{Q\in \textsf{Stop}_\beta(S)}Q\right)
\lesssim \sum_{Q\in \textsf{Stop}(S)}\sum_{Q'\in \mathcal{B}\cap
\textsf{Sons}(Q)}\mu(Q')
$$
and consequently by Lemma~\ref{lemma8.1} and the Carleson condition
(\ref{Carleson}),
\begin{align*}
\sum_{S\in\Delta_2:\;Q_S\subseteq
R}\mu(Q_S)&\lesssim\sum_{S\in\Delta_2:\;Q_S\subseteq R} \sum_{Q\in
\textsf{Stop}(S)}\sum_{Q'\in
\mathcal{B}\cap \textsf{Sons}(Q)}\mu(Q')\\
&\le \sum_{Q\in\mathcal{B}:\;Q_S\subseteq R} \mu(Q)\\
&\lesssim \mu(R).
\end{align*}

Let us consider the case $S\in \Delta^{'\phantom{''}}_{3}$. We will
need $\beta$-numbers defined for balls $B(x,r)$:
$$
\beta_q(x,r):=\inf_L
\left(\frac{1}{r}\int_{B(x,2r)}\left(\frac{\dist(x,L)}{r}\right)^q
d\mu(x)\right)^{1/q},\qquad 1\le q\le \infty,
$$
where the infimum is taken over all affine lines $L$.

 It is claimed in \cite[Section 12, Inequality (12.2)]{DavSem_Ast} that for all $S\in
\Delta_3$ there exists $\eta_2>1$ such that
$$
\mu(Q_S)\lesssim \iint_{X_S} \beta_1(x,\eta_2 r)^2\frac{d\mu(x)dr}{r},
$$
where
\begin{equation}
\label{X_S}
X_S:=\{(x,r)\in \spt \mu\times \mathbb{R}^+: x\in \eta_2Q_S,
\tfrac{1}{\eta_2 }d(x)\le r\le \eta_2 \diam (Q_S)\}
\end{equation}
and
$$
d(x):=\inf_{Q\in S}\{\dist(x,Q)+\diam(Q)\}.
$$

By Holder's inequality, $\beta_1(x,\eta_2 r) \lesssim \beta_2(x, \eta_2r)$.
Moreover, it follows from \cite[Lemma~2.5 and Proof of
Proposition 2.4]{L} (or by the arguments analogous to the ones in the proof of~Lemma~\ref{lemma_new}) that if $\mu$ is AD-regular, then there exists
$\eta_3\ge 2$ such that for any $x\in \spt \mu$,
$$
\beta_2(x,\eta_2r)^2\lesssim
\frac{1}{\eta_2r}\iiint_{\mathcal{O}_{\eta_3}(x,\eta_2r)}
c(u,v,w)^2d\mu(u)d\mu(v)d\mu(w),
$$
where
$$
\mathcal{O}_{\eta_3}(x,\rho):=\left\{(u,v,w)\in (B(x,\eta_3
\rho))^3:|u-v|\ge \frac{\rho}{\eta_3},|v-w|\ge \frac{\rho}{\eta_3},|u-w|\ge \frac{\rho}{\eta_3}\right\}.
$$
Note also that for any $(u,v,w)\in
\mathcal{O}_{\eta_3}(x,\rho)$ we have $|u-v|\le 2\eta_3\rho$, $|v-w|\le 2\eta_3\rho$ and
$|u-w|\le 2\eta_3\rho$, and thus for a fixed $\eta_3$,
$$
|u-v| \approx |v-w| \approx |u-w| \approx \rho.
$$
Therefore if a triple $(u,v,w)\in \mathcal{O}_{\eta_3}(x,\eta_2r)$
with $(x,r)\in X_S$, then at
least one side of the triangle $(u,v,w)$ makes a big angle with the
vertical and horizontal. Indeed, by construction, if $\eta_1$ is chosen much bigger than $\eta_2$, then
$\beta_\infty(x,\eta_2 r)\lesssim \sqrt{\varepsilon}$, and consequently
the angle between one side of $(u,v,w)$ and the best approximating
line $L_{Q_S}$ is less than $C(\eta_3)\sqrt{\varepsilon}$ with some
$C(\eta_3)>0$. Furthermore,
$\theta_0\le\theta_V(L_{Q_S})\le\pi/2-\theta_0$ and thus the angle that one side of
$(u,v,w)$ makes with the vertical and horizontal belongs to
$$
\big(\tfrac{9}{10}\theta_0-C(\eta_3)\sqrt{\varepsilon};\pi/2-\tfrac{9}{10}\theta_0+C(\eta_3)\sqrt{\varepsilon}\big)\supseteq
\big(\tfrac{1}{2}\theta_0;\pi/2-\tfrac{1}{2}\theta_0\big),
$$
where $\varepsilon$ is chosen sufficiently small.
This fact enables us to use the clause (\textit{ii}) of
Lemma~\ref{lemma_ends} and exchange the curvature for our permutation $p_{K_t}$:
$$
\beta_2(x,\eta_2 r)^2\lesssim_{\theta_0}
\frac{1}{\eta_2r}\iiint_{\mathcal{O}_{\eta_3}(x,\eta_2r)}
p_{K_t}(u,v,w)\;d\mu(u)d\mu(v)d\mu(w),\qquad (x,r)\in X_S.
$$
Summarizing, we get
$$
\mu(Q_S)\lesssim \iint_{X_S}  \iiint_{\mathcal{O}_{\eta_3}(x,\eta_2r)}
p_{K_t}(u,v,w)\;d\mu(u)d\mu(v)d\mu(w)\;\frac{d\mu(x)dr}{(\eta_2r)^2}.
$$
What is more, it is shown after  \cite[Lemma 7.9]{DavSem_Ast} that
the regions $X_S$ (see (\ref{X_S})) with $S\in\Delta_3^{'}$ have bounded
overlap. By this reason,
\begin{align*}
&\sum_{S\in\Delta_3^{'}:\:Q_S\subseteq R}\mu(Q_S)\\
&\quad \lesssim
\sum_{S\in\Delta_3^{'}:\:Q_S\subseteq R}\iint_{X_S}
\iiint_{\mathcal{O}_{\eta_3}(x,\eta_2 r)}
p_{K_t}(u,v,w)\;d\mu(u)d\mu(v)d\mu(w)\;\frac{d\mu(x)dr}{(\eta_2r)^2}\\
&\quad\lesssim \int_0^{2\eta_2\ell(R)}\int_{2\eta_2R}
\iiint_{\mathcal{O}_{\eta_3}(x,\eta_2r)}
p_{K_t}(u,v,w)\; d\mu(u)d\mu(v)d\mu(w)\; \frac{d\mu(x)dr}{(\eta_2r)^2}\\
&\quad\lesssim_{\eta_2} \iiint_{(u,v,w)\in (2\eta_2 R)^3} p_{K_t}(u,v,w)\; d\mu(u)d\mu(v)d\mu(w)\\
&\quad= p_{K_t}(\mu \lfloor (2\eta_2 R)).
\end{align*}
The third inequality is by Fubini's theorem. See the definition of $\mathcal{O}_\tau$ in (\ref{O_tau}). Finally, by the
$L^2$-boundedness of $T_{K_t}$, we get
$$
\sum_{S\in\Delta_3^{'}:\:Q_S\subseteq R}\mu(Q_S) \lesssim  \mu(R).
$$

Suppose now that $S\in\Delta_3^{''}$. If $Q\in
\textsf{Stop}_\alpha(S)$, then $\textsf{Sons}(Q)\cap \mathcal{B}=\varnothing$
and by Lemma~\ref{lemma8.1},
$$
\measuredangle(L_Q,L_{Q_S})\le \alpha(S),\qquad
  \measuredangle(L_Q,L_{Q_S})\ge \tfrac{1}{2}\alpha(S),\qquad \alpha(S)=10\theta_0,
$$
and thus
\begin{align*}
&\theta_V(L_Q)\le \measuredangle(L_Q,L_{Q_S})+\theta_V(L_{Q_S})<10\theta_0+\theta_0=11\theta_0,\\
&\theta_V(L_Q)\ge
\measuredangle(L_Q,L_{Q_S})-\theta_V(L_{Q_S})>5\theta_0-\theta_0=4\theta_0.
\end{align*}
Since $\beta_1(Q)<\varepsilon$, we can choose $\varepsilon$ small
enough in order that $\measuredangle(L_Q,L_{Q'})\le
\theta_0$, $Q'\in \textsf{Sons}(Q)$, and hence
$$
3\theta_0<\theta_V(L_{Q'})<12\theta_0,\qquad Q'\in \textsf{Sons}(Q).
$$
Consequently, any element of $\textsf{Sons}(Q)$ is the maximal $\mu$-cube of a
tree belonging either to $\Delta_1$, $\Delta_2$ or $\Delta_3^{'}$.
Additionally, from the definition of $\Delta_3$ and the fact that
minimal cubes for a single tree are pairwise disjoint it follows
that
$$
\mu(Q_S)\le4 \mu\left(\medcup_{Q\in
\textsf{Stop}_\alpha(S)}Q\right)= 4 \sum_{Q\in
\textsf{Stop}_\alpha(S)}\mu(Q) =4 \sum_{Q\in
\textsf{Stop}_\alpha(S)}\sum_{Q'\in \textsf{Sons}(Q)}\mu(Q').
$$
From the above-mentioned we deduce that
$$
\sum_{S\in\Delta_3^{'}:\; Q_S\subseteq R}\mu(Q_S)\le 4
\sum_{S\in\Delta_3^{'}:\; Q_S\subseteq R\phantom{L^{l}}}\sum_{Q\in
\textsf{Stop}_\alpha(S)\phantom{L^{l}}}\sum_{Q'\in \textsf{Sons}(Q)}\mu(Q') \le 4
\sum_{S\in \Delta_1\cup \Delta_2\cup \Delta_3^{'}}\mu(Q_S).
$$
Take into account that the maximal cubes of all trees from
$\Delta_1\cup \Delta_2\cup \Delta_3^{'}$ satisfy a Carleson packing
condition. By this reason,
$$
\sum_{S\in\Delta_3^{'}:\; Q_S\subseteq R}\mu(Q_S)\lesssim \mu(R).
$$

One can argue for $S\in\Delta_3^{'''}$ in the same manner as for $S\in\Delta_3^{''}$. Indeed, if $\varepsilon$ is
 appropriately chosen and $Q\in \textsf{Stop}_\alpha(S)$, then
$$
\pi/2 -12\theta_0<\theta_V(L_{Q'})<\pi/2 -3\theta_0, \qquad  Q'\in \textsf{Sons}(Q),
$$
and hence any element of $\textsf{Sons}(Q)$ is the maximal $\mu$-cube of a tree
belonging either to $\Delta_1$, $\Delta_2$ or $\Delta_3^{'}$.

Summarizing, we proved that maximal cubes of all types of trees
satisfy a Carleson packing condition and so the triple
$(\mathcal{B},\mathcal{G},\textsf{Tree})$ is a corona decomposition
as required.

\section{Additional remarks}

To finish, we would like to mention a corollary of the results
presented in the previous sections. Let $\mu$ be a Radon measure on
$\mathbb{C}$ with linear growth. If the CZO associated with the kernel $k_t$, where
 \begin{equation}
 \label{range_t_remark}
 t\in(-\infty;-\sqrt{2})\cup
(0;\infty),
\end{equation}
is $L^2(\mu)$-bounded, then all $1$-dimensional
CZOs associated with odd and sufficiently smooth kernels are also
$L^2(\mu)$-bounded. We refer the reader to \cite[Sections
1~and~12]{Tolsa2004} for the more precise description of what we mean by ``sufficiently smooth kernels''.

Indeed, it follows from~(\ref{norms}) and (\ref{norms_K_t}) with
$(n,N)=(1,2)$ that for any $t$ as in (\ref{range_t_remark}) and any cube $Q\subset \mathbb{C}$, one has
$$
\|T_{\kappa_1,\varepsilon}
\chi_Q\|_{L^2(\mu\lfloor Q)}\le C(t)\left(\|T_{k_t,\varepsilon}\chi_Q\|_{L^2(\mu\lfloor Q)}+ \sqrt{\mu(Q)} \right),\qquad C(t)>0,
$$
where $T_{\kappa_1}$, as we have already mentioned before, is the
CZO associated with the real part of the Cauchy kernel, i.e. with
the Cauchy kernel, up to a constant.  Using the $T1$ Theorem from
\cite[Theorem 9.42]{Tolsa_book}, we conclude that the $L^2(\mu)$-boundedness of $T_{k_t}$
with $t$ as in~(\ref{range_t_remark}) implies that the Cauchy
transform is $L^2(\mu)$-bounded. Furthermore, as proved in \cite{Tolsa2004}, if the
Cauchy transform is $L^2(\mu)$-bounded, then all $1$-dimensional
CZOs associated with odd and sufficiently smooth kernels are also
$L^2(\mu)$-bounded.

\end{document}